# Unified Music Theories for General Equal-Temperament Systems

Brandon Tingyeh Wu

Research Assistant,
Research Center for Information Technology Innovation, Academia Sinica,
Taipei, Taiwan

**ABSTRACT**

Why are white and black piano keys in an octave arranged as they are today? This article examines the relations between abstract algebra and key signature, scales, degrees, and keyboard configurations in general equal-temperament systems. Without confining the study to the twelve-tone equal-temperament (12-TET) system, we propose a set of basic axioms based on musical observations. The axioms may lead to scales that are reasonable both mathematically and musically in any equal-temperament system. We reexamine the mathematical understandings and interpretations of ideas in classical music theory, such as the circle of fifths, enharmonic equivalent, degrees such as the dominant and the subdominant, and the leading tone, and endow them with meaning outside of the 12-TET system. In the process of deriving scales, we create various kinds of sequences to describe facts in music theory, and we name these sequences systematically and unambiguously with the aim to facilitate future research.

## 1. INTRODUCTION

**Keyboard configuration and combinatorics**

The concept of key signatures is based on keyboard-like instruments, such as the piano. If all twelve keys in an octave were white, accidentals and key signatures would be meaningless. Therefore, the arrangement of black and white keys is of crucial importance, and keyboard configuration directly affects scales, degrees, key signatures, and even music theory. To debate the key configuration of the twelve-tone equal-temperament (12-TET) system is of little value because the piano keyboard arrangement is considered the foundation of almost all classical music theories. Modes such as Ionian and Aeolian, degrees such as the dominant and the mediant, and chords such as the major triad and the seventh chord are all based on this foundation. Beyond the 12-TET system, however, the field is not unified. Due to the rareness of pieces written in non-12-TET systems, and the fact that most non-12-TET music is generated on computers, the keyboard layout and the naming system for general equal-temperament systems remain unresolved. The starting point of this article is therefore the following question: "If there were a piano in *n* equal temperament, what should it look like?"

We can derive some basic rules to start our search of possible key combinations. Firstly, we can set the first key in an octave to be a white key without losing generality. In addition, the arrangement of white and black keys (black on white) prohibits adjacent black keys. This implies that the number of white keys in an octave can never be less than that of black keys.

Let $N_n$ be the number of possible combinations restricted by the rules above for an octave with *n* keys.

For *n* = 1, a single white key is the only solution, denoted as {○}.
For *n* = 2, a set of solutions {○○, ○●} exists.
For *n* = 3, {○○○, ○●○, ○○●}.
For *n* = 4, {○○○○, ○●○○, ○○●○, ○○○●, ○●○●}.
…
It is easy to discover that
$N_4 = C_0^4 + C_1^3 + C_2^2 = 5$, and by induction
$N_5 = C_0^5 + C_1^4 + C_2^3 = 8$, and $N_6 = C_0^6 + C_1^5 + C_2^4 + C_3^3 = 13$, and so on.

We can write $N_n = \sum_{i=0}^{\lfloor \frac{n}{2} \rfloor} C_i^{n-i}$, resulting in a Fibonacci sequence.

Specifically, $N_n = F_{n+1}$, where $F_1 = 1$, $F_2 = 1, F_n = F_{n-1} + F_{n-2}$.
The proof is straightforward:

For $N_{2k} + N_{2k+1}, k = 0,1,2,…$

$N_{2k} + N_{2k+1} = (C_0^{2k} + C_1^{2k-1} + \cdots + C_k^k) + (C_0^{2k+1} + C_1^{2k} + \cdots + C_k^{k+1}) = C_0^{2k+1} + (C_0^{2k} + C_1^{2k}) + \cdots + (C_{k-1}^{k+1} + C_k^{k+1}) + C_k^k = 1 + C_1^{2k+1} + \cdots + C_k^{k+2} + 1 = C_0^{2k+2} + C_1^{2k+1} + \cdots + C_k^{k+2} + C_{k+1}^{k+1} = N_{2k+2}.$

By the same method, we can prove that $N_{2k+1} + N_{2k+2} = N_{2k+3}$, which results in a Fibonacci sequence.

Assume that we are interested in $N_{12} = 233$. The solutions range from ○○○○○○○○○○○○ to ○●○●○●○●○●○●, and of course the regular arrangement of piano keys, ○●○●○○●○●○●○, is one of the solutions. The difficulty here is that the number of solutions grows exponentially with *n*, as

$N_n = F_{n+1} \sim \frac{1}{\sqrt{5}} \left(\frac{\sqrt{5}+1}{2}\right)^{n+1}.$

Clearly, our simple rules are not adequate. The challenge lies in making meaningful rules to restrict the freedom of key arrangement in a wise manner. Before we proceed with this, we must first consider a key part of music theory, namely, the key signature.

# 2. KEY SIGNATURES AND ABSTRACT ALGEBRA

**Key and Integers**

If you ask a musician "Which major has a key signature with four sharps?" the instant answer will always be "E major," but why is that so? In music theory, the relationship between major and minor keys, and their corresponding key signatures, can be explained by the circle of fifths. In short, step intervals of perfect fifths (or seven semitones) starting from C (i.e., C-G-D-A-E-B-…,) may lead to C major, G major, D major, etc., which have zero sharps, one sharp, two sharps, and so on. Only the "fifth" or the number "7" remain enigmatic. As will be shown later, the "7" stems from the fact that

$$n_w^{-1} = 7^{-1} \equiv 7 \ (mod \ 12)$$

where $n_w$ is the number of white keys in an octave, namely, seven.

To gain further understanding of the mathematics behind key signatures, we adopt the integer notation, which is the translation of pitch classes or interval classes into groups of integers. For 12-TET, notes are represented by integers between 0 and 11. For example, C = 0, C# = 1, G = 7, and B = 11. The integer notation is especially useful for equal-temperament systems because the intervals between these numbers have a direct connection to our perception of pitch interval. Consider the major thirds C-E and D-F#. Their intervals are both 4 since 4 – 0 = 6 – 2 = 4. So far, a set of numbers and an operation "+" ("+" or "-" are equivalent ideas) are involved. The use of the operation "+" needs to be clarified because, for instance, 7 + 7 = 14, but 14 is not in our integer set. When using "+," we are actually computing the note that occurs at a certain interval above another note. The perfect fifth above G = 7 is D = 2, so we want 7 + 7 to be 2. In this sense, the algebraic structure that we adopt is the cyclical-quotient group (or more precisely, "ring of integers" when operation of multiplication modulo 12 become available later) with the operation-addition modulo 12, or in mathematical terms, $\mathbf{Z}/12\mathbf{Z}$, or simply $\mathbf{Z}_{12}$.

**Computation of Key Signatures**

To compute the key signatures, we can follow a simple procedure. First, we need the positions of the white keys, or the C major diatonic set. Clough and Douthett (1991) provided a clear definition, and notated a subset of $d$ pitch classes selected from the chromatic universe of $c$ pitch classes as $\mathrm{D}_{c,d}$. Following this notation, in classical 12-TET, $\mathrm{D}_{12,7} = \{0,2,4,5,7,9,11\}$. However, in our article, transposition and sequence sorting are widely used; hence, we define a sequence as $S_j^{(n,n_w),i} = \{j\} + \mathrm{D}_{n,n_w} \equiv \{j + D_1, j + D_2, \ldots, j + D_{n_w}\} \ (mod \ n)$, with elements in ascending order.

The $i$ in the superscript is an indicator of specific sequences among possible $S_j^{(n,n_w)}$. As will be explained in parts III and IV, there are three $S_j^{(n,n_w)}$ for $(n, n_w) = (12,7)$, and $D_{12,7} = \{0,2,4,5,7,9,11\}$ equals $S_0^{(12,7),2}$. We omit $i$ here, and if $S_j$ is not specified, it refers to $S_j^{(12,7),2}$. To find a key signature, for example, B major, we first find its tonic (e.g., B = 11), and add it to the sequence as $11 + \{0,2,4,5,7,9,11\} \equiv \{11,1,3,4,6,8,10\}\ (mod\ 12)$. We then sort the sequence to get $S_{11} = \{1,3,4,6,8,10,11\}$, and subtract the C major scale sequence, resulting in the key signature sequence $K_{11} = S_{11} - S_0 = \{1,3,4,6,8,10,11\} - \{0,2,4,5,7,9,11\} = \{+1,+1,0,+1,+1,+1,0\}$. This corresponds to a key with five sharps, namely C#, D#, F#, G#, and A#, which is indeed the B-major scale. In the following example, we compute the A-flat major scale as

$$8 + \{0,2,4,5,7,9,11\} \equiv \{8,10,0,1,3,5,7\} \ (mod\ 12)$$

and $K_8 = S_8 - S_0 = \{0,1,3,5,7,8,10\} - \{0,2,4,5,7,9,11\} = \{0,-1,-1,0,0,-1,-1\}$. This corresponds to a key with four flats, namely, Db, Eb, Ab, and Bb, which is indeed the A-flat scale.

We do not rewrite the negative value (i.e., "-1") as "11" for several reasons, despite the fact that "-1" is not contained in our $Z_{12}$ group.

Firstly, since we use "+1" or simply "1" to denote the location of sharps, it is straightforward to denote flats as "-1."

Secondly, if we define the function $\|K_n\| = sum\ of\ elements\ in\ K_n$, then a positive value gives the number of sharps in the key signature, and similarly a negative value gives the number of flats, such as $\|K_{11}\| = +5$ and $\|K_8\| = -4$ (to avoid confusion with the use of the set function "cardinality," we use the notation of $\|K_n\|$ instead of $|K_n|$). Another reason for not rewriting negative values is closely connected to the issue of enharmonic equivalent. In some cases, we want -1 and 12 to be the aliases of 11 and 0, respectively.

**Enharmonic Equivalence**

If we consider the commonly seen major key signature, from seven sharps to seven flats, then there should be 7+1+7 = 15 different keys. However, there are only 12 different major scales from the perspective of piano performance. The discrepancy stems from the fact that B major is enharmonically equivalent to C-flat major; F-sharp major is enharmonically equivalent to G-flat major; and C-sharp major is enharmonically equivalent to D-flat major. In other words, it is possible that all notes on two scores are different, but if this is accompanied by a change in key signature, no difference can be perceived during performance. Examples of this can be found in some of Franz Liszt's works. For example, "La Campanella," as the third of Franz Liszt's six Grandes etudes de Paganini ("Grand Paganini Etudes," opus S. 141, 1851) is called, was written in G-sharp minor, but its older version, Etudes d'execution transcendente d'apres Paganini (opus S. 140, 1838), was written in A-flat minor.

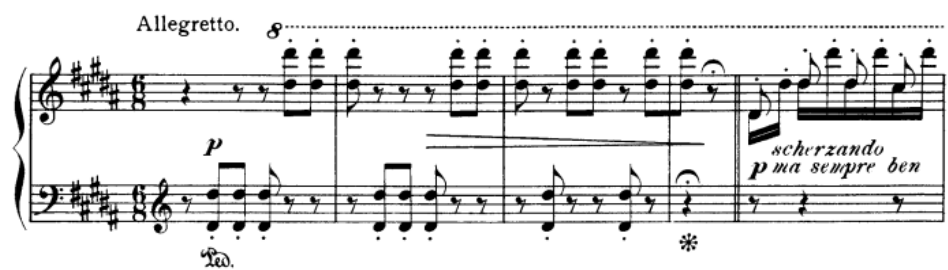

Example 1a. Liszt's Grandes études de Paganini, S.141

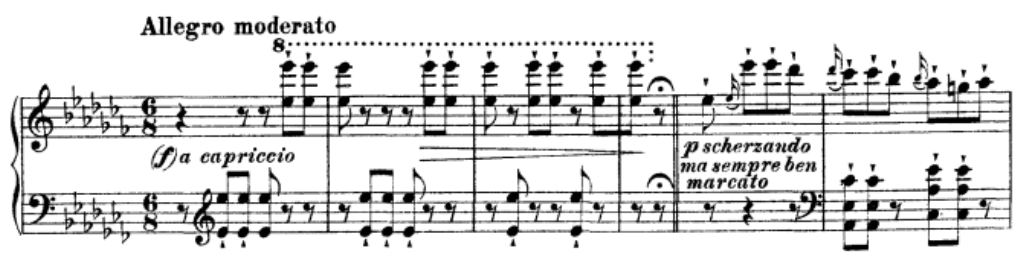

Example 1b. Liszt's Études d'exécution transcendante d'après Paganini, S.140

Another example is Liszt's Hungarian Rhapsody No. 6, where the first part (mm. 1-40) was written in D-flat major, but the presto part (mm. 41-72) was written in C-sharp major. However, no key change is perceived by the human ear.

Example 2. Liszt's Hungarian Rhapsody No. 6, S.244/6

These two examples serve to illustrate the enharmonic equivalence of two different key signatures. To make our key signature sequence $K_n$ capable of indicating enharmonic equivalence, we allow -1 and 12 to be the aliases of 11 and 0, respectively, in $Z_{12}$.

Accordingly, before we sort the sequence $\{11,1,3,4,6,8,10\}$ in B major, we can take -1 as the alias of 11, having $S_{11}^{(-)} = \{-1,1,3,4,6,8,10\}$. Then $K_{11}^{(-)} = S_{11}^{(-)} - S_0 = \{-1,1,3,4,6,8,10\} - \{0,2,4,5,7,9,11\} = \{-1,-1,-1,-1,-1,-1,-1\}$, resulting in a key signature with seven flats, which is enharmonically equivalent to a key signature with five sharps. The relation between two enharmonically equivalent keys is closely connected to the fact that $\|K_{11}\| = +5 \equiv \|K_{11}^{(-)}\| = -7 \equiv 5 \ (mod\ 12)$.

When considering the example of C-sharp major, we can take 12 as the alias of 0, having $K_1^{(+)} = S_1^{(+)} - S_0 = \{1,3,5,6,8,10,12\} - \{0,2,4,5,7,9,11\} =$ $\{+1, +1, +1, +1, +1, +1, +1\}$. This results in a key signature with seven sharps, which is enharmonically equivalent to a key signature with five flats, and of course $\|K_1\| = -5 \equiv \|K_1^{(+)}\| = +7 \equiv 7 \ (mod\ 12)$.

However, including -1 and 12 in the sequence may result in some theoretical keys that are rarely seen in practice. Consider, for example, $8 + \{0,2,4,5,7,9,11\} \equiv \{8,10,0,1,3,5,7\}$ in A-flat major. If we allow 12 to be the alias of 0 and re-sort the sequence, we have $S_8^{(+)} = \{1,3,5,7,8,10,12\}$. When we compare this to $S_0$, we have $K_8^{(+)} = S_8^{(+)} - S_0 = \{1,3,5,7,8,10,12\} - \{0,2,4,5,7,9,11\} = \{+1, +1, +1, +2, +1, +1, +1\}$ and $\|K_8^{(+)}\| = +8$. The +2 in $K_8^{(+)}$ can be interpreted as a double sharp, and the resulting key is G-sharp major, the enharmonic equivalence of A-flat major. Scores with G-sharp major can be found in John Foulds' "A World Requiem" (1921).

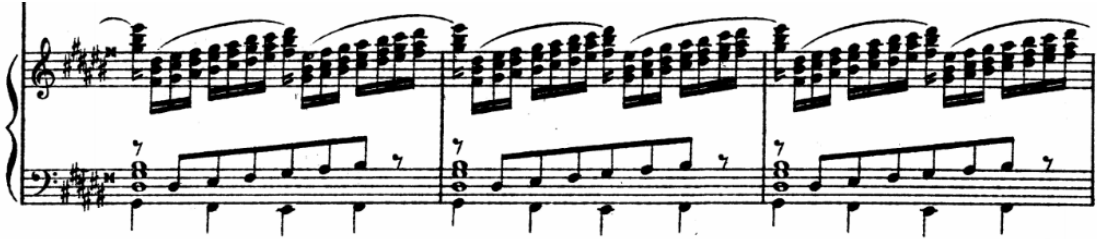

Example 3. John Foulds' "A World Requiem, Op. 60"

**Generating Sequences**

As explored by Clough and Myerson (1985) and defined by Clough and Douthett (1991), a $k$-generated maximally even (ME) set, $M_{n,n_w}^k$, is a set that can be represented as $\{j, j+k, j+2k, \ldots, j+(n_w-1)k\}, 0 \leq j, k \leq n-1$. Based on this idea, we define a sequence as $G_{(s_0,k),m}^{(n,n_w)} \equiv \{s_0, s_0+k, \ldots, s_0+(m-1)k\} \ (mod\ n)$. Note that $m$ is the length of the sequence and is not necessarily $n_w$.

$S_0$ is a subset of $\mathbf{Z}_{12}$ and is an ME set that can be generated with $(s_0, k) = (5,7)$ and $(11,5)$. Determining the number of generators for a scale is not a trivial question, and a careful survey has been conducted by Amiot (2009). However, there are always at least two generators, because whenever $k$ is a generator, $n-k$ is also a generator. The generator pair $(s_0, k) = (5,7)$ and $(11,5)$ corresponds to the generating sequence $\{5,0,7,2,9,4,11\}$ and $\{11,4,9,2,7,0,5\}$, respectively. Note that when $mk \equiv 1 \ (mod\ n)$, we can extend our sequence with $m$ elements both forwards, $G_{(s_0,k),m}^{(n,n_w)\ (+)} \equiv G_{(s_0+mk,k),m}^{(n,n_w)} = G_{(s_0+1,k),m}^{(n,n_w)}$, and backwards, $G_{(s_0,k),m}^{(n,n_w)\ (-)} \equiv G_{(s_0-mk,k),m}^{(n,n_w)} = G_{(s_0-1,k),m}^{(n,n_w)}$, to get $G_{(s_0-1,k),3m}^{(n,n_w)}$. To solve $mk \equiv 1 \ (mod\ n)$ for a given $m$ is to find the modular multiplicative inverse $k$ of $m$ in the ring of integers modulo $n$, and $k$ can be expressed as $k \equiv m^{-1} \ (mod\ n)$. $k$ exists if and only if $m$ and $n$ are coprime, and thus $k$ and $n$ will be coprime as well. In the 12-TET case, among all possible pairs of generators $(s_0, k)$, only $(s_0, k) = (5,7)$ satisfies $mk \equiv 1 \ (mod\ n)$, namely, $7 \cdot 7 \equiv 1 \ (mod\ 12)$. Thus, $G_{(s_0-1,k),3m}^{(n,n_w)} = G_{(4,7),21}^{(12,7)} =$ $\{4,11,6,1,8,3,10,5,0,7,2,9,4,11,6,1,8,3,10,5,0\}$. In the following, $(n, n_w)$ may be omitted as they refer to $(12,7)$. The purpose of determining $G_{(s_0-1,k),3m}$ is that if $G_{(s_0,k),m}$ contains exactly all $m$ elements in S, then $G_{(s_0+k,k),m}$ contains exactly all $m$ elements in $k + S$. Recall that $S_0 = \{0,2,4,5,7,9,11\}$ is contained in $G_{(5,7),7} = \{5,0,7,2,9,4,11\}$, and $S_7 = \{0,2,4,6,7,9,11\}$ can be contained in $G_{(5+7,7),7} = G_{(0,7),7} = \{0,7,2,9,4,11,6\}$. There is an overlap of $m - 1$ elements between $G_{(s_0,k),m}$ and $G_{(s_0+k,k),m}$, and we can interpret this as a right-shift of a window of size $m$ in the sequence $G_{(s_0-1,k),3m}$, as illustrated in Figure 4.

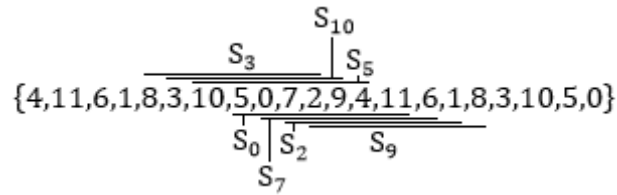

Figure 4. $G_{(s_0-1,k),3m}$ and the sequence that it contains

Note that whenever we shift the window to the right from $G_{(s,k),m}$ to $G_{(s+k,k),m}$, the elements that it meets change exactly from $S_{s'}$ to $S_{s'+k}$. In the case above, it is $s' = s - s_0 = s + 7$. $S_{s'+k}$ differs from $S_{s'}$ in only one element because the first element in $G_{(s,k),m}$ is replaced by the last element in $G_{(s+k,k),m}$, that is, $s$ is replaced by $s + mk = s + 1$. This is an important property in that whenever we leap from $S_{s'}$ to $S_{s'+k}$, the key signature gains one sharp (or loses one flat) and we can therefore make sure no key signature is left behind. In our notation, this means $\|K_{s'+k}\| = \|K_{s'}\| + 1$, which corresponds to the "circle of fifths." $k = 7$ is the interval of a perfect fifth, and with the tonic raised by a fifth, the corresponding key signature has an additional sharp or is equivalently one flat less. Since $\|K_0\| = 0$ by $K_0 = S_0 - S_0 = \{0,0,0,0,0,0,0\}$, we have $\|K_k\| = \|K_{m^{-1}}\| = 1$, $\|K_{2k}\| = \|K_{2(m^{-1})}\| = 2, \dots$, $\|K_{mk}\| = \|K_{m \cdot m^{-1}}\| = \|K_1\| = m$, and $\|K_t\| = tm$, which gives the relation between tonics and key signatures. The position of sharps or flats can then be induced by the sequences $G_{(s_0,k),m}^{(+)}$ and $G_{(s_0,k),m}^{(-)}$ (reverse), respectively.

In order to eliminate ambiguity, $(s_0, k)$ is chosen uniquely so that $mk \equiv$

$1 \ (mod \ n)$. In addition, the $m$ of interest is the number of white keys in an octave, $n_w$. So for a given scale with $n$ keys and $n_w$ white keys, $G_{(s_0,k),m} = G_{(s_0,n_w^{-1}),n_w}$ if exists, is unique. Since the number $n_w^{-1}$ has a special position in any $n$-TET systems, we can define scale degrees as follows:

**Tonic:**

  The tonic is the main or first note of the scale, so the tonic is the note with number 0.

**Dominant:**

  The conventional definition of the dominant is that it is the perfect fifth of the tonic (i.e., note number 7). We can broaden this definition to make the note with number $n_w^{-1}$ the dominant. By this definition, the dominant key is not always the key whose tonic is a perfect fifth (+7) above the tonic of the main key, but $n_w^{-1}$ above the tonic of the main key. Despite this fact, $n_w^{-1}$ still has a very strong connection with the perfect fifth in some systems, such as 19-TET and 31-TET (for more details, please see the universal perfect-fifth constant in part VI). This guarantees that the signature of the dominant key has always one more sharp or one less flat than the original key. In the traditional 12-TET system, this definition of the dominant coincides with the convention since $n_w^{-1} = 7 \ (mod \ 12)$.

**Subdominant:**

  Similar to the definition of the dominant, we define the subdominant as the note with number $-n_w^{-1}$, or equivalently $n - n_w^{-1}$. In doing so, the signature of the subdominant key always has one less sharp or one more flat than the original key.

**Leading tone:**

  In the conventional 12-TET system, $G_{(5,7),7} = \{5,0,7,2,9,4,11\}$. If we add 7 (and of course mod 12) to any of its elements except the last, the new number is still in the sequence $S_0 = \{0,2,4,5,7,9,11\}$. However, if we add 7 to the last element of 11, we get 6, which is not included in $S_0$. This means that the last element in $G_{(5,7),7}$ "leads" the sequence $G_{(5,7),7}$ to "leave" $S_0$ by having the first sharp. Indeed, 6 (F#) is the first sharp we see in the key signature. We can therefore define the last element in $G$ as the leading tone. This definition differs from the traditional one, which states that the leading note is the seventh tone or degree of a scale that is a half tone below the tonic.

  Let us apply this argument to the first element in $G_{(5,7),7}$. If we subtract 7 from any but the first of the elements in $G_{(5,7),7}$, the new number will still be in the sequence $S_0 = \{0,2,4,5,7,9,11\}$. So the first element in $G_{(5,7),7}$ leads the sequence $G_{(5,7),7}$ to leave $S_0$ if we take a descending step. Therefore, we can categorize leading tones into either ascending or descending leading tones, which are defined as follows:

**Descending leading tone:**
　　The note with number $s_0$, which is the first element in $G_{(s_0,n_w^{-1}),n_w}$.

**Ascending leading tone (conventional leading tone):**
　　The note with number $s_0 - n_w^{-1} + 1$, which is the last element in $G_{(s_0,n_w^{-1}),n_w}$.

　　Please note that the ascending leading tone in our definition coincides with the conventional leading tone in 12-TET since $s_0 - n_w^{-1} + 1 = 5 - 7 + 1 = -1 \equiv 11 \pmod{12}$.

# 3. KEY CONFIGURATIONS

Based on the above observations and arithmetic background on key signatures, we now propose some axioms based on established musical facts, and determine the key arrangements in other TET systems.

**Axioms**

As described in part I, if we do not limit the ratio of black and white keys, the number of possible scales grows drastically with the number of keys in an octave. The first axiom serves to regulate the relation between neighboring keys. Assuming that we do not want our diatonic scale to be too "chromatic," i.e., having many consecutive white keys not separated by black keys, we can define our first axiom for the key arrangement in TET systems as follows.

**Axiom I:**

**Three consecutive white keys not separated by black keys are prohibited. In addition, there must be at least one white key between any two black keys.**

The first part of this axiom limits the lower bound of $n_b$, the number of black keys in an octave, to be $\left\lceil \frac{1}{3}n \right\rceil$, while the second part sets the upper bound of $n_b$ to be $\left\lfloor \frac{1}{2}n \right\rfloor$. The second part makes the chromatic length of a step either one or two, resulting in a reduced and rounded scale, as defined by Clough and Myerson (1986). Although this axiom does not alter the exponential nature of the number of possible solutions, it eliminates the number of possible key arrangements considerably, so that the number of solutions does not go out of control before *n* grows too large for people to discern the pitch of two adjacent keys.

The mathematical representations developed in part II may help us to define the next axiom. From $K_0 = \{0,0,0,0,0,0,0\}$, $K_1 = \{0,-1,-1,0,-1,-1,-1\}$ to $K_{11} = \{+1,+1,0,+1,+1,+1,0\}$, we can see that all $K_i$ are spanned by either $\{0,1\}$ or $\{0,-1\}$, and that $K_0$ is a special case spanned by $\{0\}$. When $K_i$ is spanned by $\{0\}$, we call $K_i$ monotonic, and when $K_i$ is spanned by $\{0,1\}$ or $\{0,-1\}$, we call $K_i$ a monotonic single sharp or flat sequence, respectively. Our second axiom addresses the monotonicity of $K_i$.

**Axiom II:**

**For a given number of keys in an octave ($n$), a permitted key arrangement must be such that $K_i$ is either a monotonic single sharp sequence or a monotonic single flat sequence, for $i = 1, \dots, n-1$.**

This axiom states that the key signature of any given key should contain either only sharps or only flats. Conventional music scores never have key signatures that contain both sharps and flats, and this is indeed the starting point of this axiom. The following is an example of the key arrangement that violates Axiom II but not Axiom I:

For n = 12, consider {○●○○●○●○●○●○}, that is, $S_0 = \{0,2,3,5,7,9,11\}$, which is the ascending melodic minor scale in the original 12-TET system. $K_5 = S_5 - S_0 = \{0,2,4,5,7,8,10\} - \{0,2,3,5,7,9,11\} = \{0,0,+1,0,0,-1,-1\}$, and it is neither a monotonic single sharp nor a monotonic single flat.

It follows that the property of maximal evenness introduced by Clough and Douthett (1991) is a necessary condition of this axiom; otherwise, there must be +2 or -2 in some $K_i$. However, the stricter Myhill's Property of symmetry is not implied by this axiom.

A further axiom is required to address the problem of enharmonic equivalence, in which a key may be represented by different key signatures. To exclude cases in which different keys map to the same key signature, we define the following axiom.

**Axiom III:**

**For $i - j \neq 0 \ (mod \ n)$, $K_i - K_j \neq K_0$, or equivalently, the map from tonic $t$ to key signature $K_t$ is injective.**

When we see a key signature, we are able to specify the key and its tonic, the dominant, and the leading tone without ambiguity. When finding key arrangements outside of the 12-TET system, we want the property of injective mapping to hold true. Usually, only a few scales abide by Axioms I and II but violate Axiom III, but the periodic structure within an octave may violate Axiom III. Consider the following example:

For n = 15, consider {○●○●○○●○●○○●○●○}, that is, $S_0 = \{0,2,4,5,7,9,10,12,14\}$, so that we have $K_0 = K_5 = K_{10}$, $K_1 = K_6 = K_{11}, \dots, K_4 = K_9 = K_{14}$. In reality, there are only five distinct key signatures. In the mathematical sense, this behavior of the inner period violates the coprime property; thus, the inverse element of $n_w$ does not exist. This is a good example demonstrating that Axioms I and II do not imply Myhill's Property, since a diatonic length of 3 implies a chromatic length of 5 in this case.

# 4. KEY CONFIGURATION FOR SOME CASES

In this section, we first list the number of permitted key layouts $N_{n,n_w}$ based on our axioms for some *n*, and then list their properties.

Table 5. Number of solutions for some *n*

| $n$ | $n_w$ | $n_b$ | $N_{n,n_w}$ | Cyclically equivalent? |
|---|---|---|---|---|
| 3 | 2 | 1 | 2 | Yes |
| 5 | 3 | 2 | 2 | Yes |
| 7 | 4 | 3 | 2 | Yes |
| 8 | 5 | 3 | 3 | Yes |
| 9 | 5 | 4 | 2 | Yes |
| 11 | 7 | 4 | 4 | Yes |
| 11 | 6 | 5 | 2 | Yes |
| 12 | 7 | 5 | 3 | Yes |
| 13 | 8 | 5 | 4 | Yes |
| 13 | 7 | 6 | 2 | Yes |
| 14 | 9 | 5 | 5 | Yes |
| 15 | 8 | 7 | 2 | Yes |
| 16 | 9 | 7 | 3 | Yes |
| 17 | 11 | 6 | 6 | Yes |
| 17 | 10 | 7 | 4 | Yes |
| 17 | 9 | 8 | 2 | Yes |
| 18 | 11 | 7 | 5 | Yes |
| 19 | 12 | 7 | 6 | Yes |
| 19 | 11 | 8 | 4 | Yes |
| 19 | 10 | 9 | 2 | Yes |
| 20 | 13 | 7 | 7 | Yes |
| 20 | 11 | 9 | 3 | Yes |
| 21 | 13 | 8 | 6 | Yes |
| 21 | 11 | 10 | 2 | Yes |
| 22 | 13 | 9 | 5 | Yes |
| 23 | 15 | 8 | 8 | Yes |
| 23 | 14 | 9 | 6 | Yes |
| 23 | 13 | 10 | 4 | Yes |
| 23 | 12 | 11 | 2 | Yes |
| 24 | 13 | 11 | 3 | Yes |

Two sequences $S_a$ and $S_b$ of the same length $m$ can be called cyclic equivalents if and only if there exists an offset $d$, such that for all $0 \leq i < m$, $S_a[i] = S_b[(i+d) \mod m]$. Table 5 shows that the solutions, given a pair of $(n_w, n_b)$, are all cyclic equivalents. This fact is strongly connected to properties of rounded, reduced, and maximal evenness. In the key configurations in Table 6, 0 represents a white key and 1 represents a black key. The ascending leading tone is marked with a right arrow, and the descending with a left arrow. The dominant is underlined and the subdominant is double-underlined.

Table 6. Key configurations for some interesting cases

| $n, n_w, n_b$ | Key Configuration | $S_0$ | $(s_0, n_w^{-1})$ |
|---|---|---|---|
| 12,7,5 | 0̄10101̱0⃗01010 | {0,2,4,6,7,9,11} | (0,7) |
| | 01010⃖0̱10101⃗0 | {0,2,4,5,7,9,11} | (5,7) |
| | 01010⃗0̱101001̄ | {0,2,4,5,7,9,10} | (10,7) |
| 17,11,6 | 0̄101⃗0010010010̱010 | {0,2,4,5,7,8,10,11,13,14,16} | (0,14) |
| | 0100⃖101⃗0010010̱010 | {0,2,3,5,7,8,10,11,13,14,16} | (3,14) |
| | 0100100⃖101⃗0010̱010 | {0,2,3,5,6,8,10,11,13,14,16} | (6,14) |
| | 0100100100⃖101⃗0̱010 | {0,2,3,5,6,8,9,11,13,14,16} | (9,14) |
| | 0100100100100⃖1̱01⃗0 | {0,2,3,5,6,8,9,11,12,14,16} | (12,14) |
| | 010⃗0100100100100̱01̄ | {0,2,3,5,6,8,9,11,12,14,15} | (15,14) |
| 17,10,7 | 0̄10101⃗00101001̱010 | {0,2,4,6,7,9,11,12,14,16} | (0,12) |
| | 01010⃖010101⃗001̱010 | {0,2,4,5,7,9,11,12,14,16} | (5,12) |
| | 010100101001⃖01̱01⃗0 | {0,2,4,5,7,9,10,12,14,16} | (10,12) |
| | 0101⃗0010100101001̱̄ | {0,2,4,5,7,9,10,12,14,15} | (15,12) |
| 17,9,8 | 0̄101010101010101⃗0 | {0,2,4,6,8,10,12,14,16} | (0,2) |
| | 010⃗1010101010101̱̄01 | {0,2,4,6,8,10,12,14,15} | (15,2) |
| 19,12,7 | 0̄10100100101⃖0010010 | {0,2,4,5,7,8,10,12,13,15,16,18} | (0,8) |
| | 0101⃗0010010̱01⃖010010 | {0,2,4,5,7,8,10,11,13,15,16,18} | (11,8) |
| | 0100⃖101001001010⃗010 | {0,2,3,5,7,8,10,11,13,15,16,18} | (3,8) |
| | 01001010⃗0100100⃖1010 | {0,2,3,5,7,8,10,11,13,14,16,18} | (14,8) |
| | 010010⃖01010⃗01001010 | {0,2,3,5,6,8,10,11,13,14,16,18} | (6,8) |
| | 010010010⃗01001001⃖01 | {0,2,3,5,6,8,10,11,13,14,16,17} | (17,8) |
| 19,11,8 | 0̄1010010101⃖01̱001010 | {0,2,4,6,7,9,11,13,14,16,18} | (0,7) |
| | 010101⃗001010⃖01̱01010 | {0,2,4,6,7,9,11,12,14,16,18} | (12,7) |
| | 01010⃖01010100̱101⃗010 | {0,2,4,5,7,9,11,12,14,16,18} | (5,7) |
| | 010100101⃗0̱1001010⃖01 | {0,2,4,5,7,9,11,12,14,16,17} | (17,7) |
| 19,10,9 | 0̄10101010101010101⃗0 | {0,2,4,6,8,10,12,14,16,18} | (0,2) |
| | 01̱010101010101010⃗1̄01 | {0,2,4,6,8,10,12,14,16,17} | (17,2) |

- 16 -

Table 6 shows that solutions with the same pair of $(n_w, n_b)$ are those generated with the same $n_w^{-1}$, which makes those solutions cyclic equivalents to each another.

Generally speaking, the number of possible pairs of $n_w, n_b$ is significantly lower for even $n$ because half the integers between $\left\lceil \frac{1}{3}n \right\rceil$ and $\left\lfloor \frac{1}{2}n \right\rfloor$ are even and impossible to be coprime to $n$. For $n = 4, 6, 10$, there is no possible key arrangement according to the three axioms. However, if $n$ is a prime number, all integers between $\left\lceil \frac{1}{3}n \right\rceil$ and $\left\lfloor \frac{1}{2}n \right\rfloor$ are coprime to $n$ and produce abundant possible configurations.

Of greatest interest is definitely $n = 12$ since nearly all music has been written in the 12-TET system. It is surprising that the number of possible configurations has been reduced from the $N_{12} = 233$ described in part I to $N_{12,7} = 3$ through only three axioms. Moreover, the three solutions are actually equivalent and coincide with the traditional 12-TET system. This means that the three axioms based on the observation and generalization of key signature make sense in the context of music theory.

**Nomenclature**

From Table 5, we see that $N_{n,n_w}$ is $\leq 2$ if it does not equal 0. In fact, as shown below, $N_{n,n_w} = n_w - n_b + 1$ and is thus impossible to be equal to 1. We therefore need to order the possible key configurations properly. Ordering sets is an old problem in music theory. For instance, Forte (1973) orders pitch-class sets by Forte number in a "left-packed" way. Set class 4–2 [0124] is listed before set class 4–3 [0134] because the 1's (note indicators) in "11101" are more left-packed than in "11011." Another plausible nomenclature is to regard the binary expressions of key configurations as integers, and order them by their value. For example, the sequences 010101001010, 010100101010, and 010100101001 correspond to the decimal values of 1354, 1322, and 1321, respectively. Accordingly, we can specify their corresponding $S_0$ as $S_0^{(12,7),3} = \{0,2,4,6,7,9,11\}$, $S_0^{(12,7),2} = \{0,2,4,5,7,9,11\}$, and $S_0^{(12,7),1} = \{0,2,4,5,7,9,10\}$, respectively, to avoid making $S_0$ confusing. The two methods yield the same result because binary ordering is essentially left-packing the 0s, and in the key configurations above, 0s rather than 1s are the indicators of notes in $S_0$. Ordering results in all possible $S_0$ for a given set of $n, n_w$ to range from $S_0^{(n,n_w),1}$ to $S_0^{(n,n_w),n_w-n_b+1}$. Without proving $N_{n,n_w} = n_w - n_b + 1$ formally, we prefer to show it in a more intuitive manner with the concept of key signatures and enharmonic equivalence.

**Key signature matrices**

Recall the idea of key signature sequence introduced in part II. If we keep the expression of $K_i$ and do not allow $K_i^{(+)}$ and $K_i^{(-)}$, then we can produce a key-signature matrix with size $n_w$-by-$n$ by assigning vector $K_i$ to column $i$ of the matrix. For example, the key-signature matrix for the conventional 12-TET system $\mathbf{K}^{(12,7),2}$ of scale $S_0^{(12,7),2} = \{0,2,4,5,7,9,11\}$ is

$$\mathbf{K}^{(12,7),2} = \begin{pmatrix} 0 & 0 & 1 & 0 & 1 & 0 & 1 & 0 & 0 & 1 & 0 & 1 \\ 0 & -1 & 0 & 0 & 1 & 0 & 1 & 0 & -1 & 0 & 0 & 1 \\ 0 & -1 & 0 & -1 & 0 & 0 & 1 & 0 & -1 & 0 & -1 & 0 \\ 0 & 0 & 1 & 0 & 1 & 0 & 1 & 1 & 0 & 1 & 0 & 1 \\ 0 & -1 & 0 & 0 & 1 & 0 & 1 & 0 & 0 & 1 & 0 & 1 \\ 0 & -1 & 0 & -1 & 0 & 0 & 1 & 0 & -1 & 0 & 0 & 1 \\ 0 & -1 & 0 & -1 & 0 & -1 & 0 & 0 & -1 & 0 & -1 & 0 \end{pmatrix}.$$

If we calculate the sum for each column, we get

$$\|\mathbf{K}^{(12,7),2}\| \equiv \sum_i \mathbf{K}_{i,j}^{(12,7),2} = (\|K_0\| \quad \|K_1\| \quad \dots \quad \|K_{11}\|)$$
$$= (0 \quad -5 \quad 2 \quad -3 \quad 4 \quad -1 \quad 6 \quad 1 \quad -4 \quad 3 \quad -2 \quad 5).$$

The minimum and maximum values of $\sum_i \mathbf{K}_{i,j}^{(12,7),2}$ are -5 and 6, respectively. Since we do not allow $K_i^{(+)}$ and $K_i^{(-)}$, the minimum value of $\sum_i \mathbf{K}_{i,j}^{(12,7),\cdot}$ cannot be -7 because this would require -1 in some $K_i$, which is impossible since we disallow $K_i^{(-)}$. The maximum value of $\sum_i \mathbf{K}_{i,j}^{(12,7),\cdot}$, however, can be 7, which comes from the sum of the second column, $\mathbf{K}_{i,2}^{(12,7),1}$. That is to say, $K_1^{(12,7),1} = S_1^{(12,7),1} - S_0^{(12,7),1} = \{1,3,5,6,8,10,11\} - \{0,2,4,5,7,9,10\} = \{1,1,1,1,1,1,1\}$. If we list all $\sum_i \mathbf{K}_{i,j}^{(12,7),\cdot}$, we have

$$\sum_i \mathbf{K}_{i,j}^{(12,7),1} = (0 \quad 7 \quad 2 \quad -3 \quad 4 \quad -1 \quad 6 \quad 1 \quad -4 \quad 3 \quad -2 \quad 5)$$
$$(min, max) = (-4,7)$$

$$\sum_i \mathbf{K}_{i,j}^{(12,7),2} = (0 \quad -5 \quad 2 \quad -3 \quad 4 \quad -1 \quad 6 \quad 1 \quad -4 \quad 3 \quad -2 \quad 5)$$
$$(min, max) = (-5,6)$$

$$\sum_i \mathbf{K}_{i,j}^{(12,7),3} = (0 \quad -5 \quad 2 \quad -3 \quad 4 \quad -1 \quad -6 \quad 1 \quad -4 \quad 3 \quad -2 \quad 5)$$
$$(min, max) = (-6,5).$$

Correspondingly, $N_{n,n_w} = n_w - n_b + 1$ can be explained by the fact that there are $n_w - n_b + 1$ possibilities between $(min, max) = (-n_w + 1, n_b)$ to $(-n_b +$ $1, n_w$), or, colloquially put, "you can put a stick of length $n$ in a box of length $2n_w$ at $2n_w - n + 1 = n_w - n_b + 1$ distinct integer positions." We call this argument "the argument of capacity." Essentially, $K_1^{(n,n_w),1}$ always consists of only 1s. This property is closely connected to the prime form of $S_0^{(n,n_w),\cdot}$.

**Relations between $S_0^{(n,n_w),1}$ and the prime forms of $S_0^{(n,n_w),\cdot}$**

It is a universal property that $S_0^{(n,n_w),1}$ is the inversion of the prime form, and $S_0^{(n,n_w),1}$ is the only sequence whose interval equals n-2, with other $S_0^{(n,n_w),\cdot}$ having intervals of n-1. Throughout this article, "prime form" refers to Forte's definition rather than that in John Rahn (1980) for its consistency with the left-packing method.

While $S_0^{(n,n_w),1}$ is always the inversion of the prime form, the prime form itself is never in the list of $S_0^{(n,n_w),\cdot}$. This is a result of left-packing when computing the prime form.

The prime form is the most packed form. This property allows full sharps in the key-signature sequence since the largest element in $S_0^{(n,n_w),prime}$ (the prime form of $S_0^{(n,n_w),\cdot}$) is $n-2$, and the largest element in $S_1^{(n,n_w),prime}$ is $n-1$. Correspondingly, $K_1^{(n,n_w),prime}$ is therefore consists entirely of sharps. Up to this point, the same seems to apply to $S_0^{(n,n_w),1}$. However, we now claim that the first three elements in the prime form of $S_0^{(n,n_w),\cdot}$ are always $\{0,1,3,\dots\}$.

Clearly, there exists $1 \leq i \leq n_w - 1$ (note that throughout this article, indices of elements within a sequence or a matrix start with 1 rather than 0) such that the $i^{th}$ element and the $(i+1)^{th}$ element in $S_0^{(n,n_w),prime}$, say $x_i$ and $x_{i+1}$, differ by 1. In short, a pair of adjacent and continuous numbers must occur at some point in the sequence. This is a direct consequence of Axiom I and the definition of the prime form. Next, we prove that the first such $i$ is always 1.

Assume the first such $i$ is $j \neq 1$. We then have $S_0^{(n,n_w),prime} = \{0,2,4,\dots,x_j, x_j + 1, x_j + 3, \dots, n-2\}$. Recall that Forte's algorithm puts a priority on making small numbers smaller. If we rotate $S_0^{(n,n_w),prime}$ to the left and subtract two, the resulting new sequence is more left-packed while preserving the same whole interval $n-2$. If we repeat this procedure $j-1$ times as shown in Figure 7, we will get a more (but not necessarily the most) left-packed sequence.

| 0 | 2 | 4 | ... | $x_j$ | $x_j+1$ | $x_j+3$ | ... | $n-2$ |
|---|---|---|---|---|---|---|---|---|
| 2 | 4 | ... | $x_j$ | $x_j+1$ | $x_j+3$ | ... | $n-2$ | $n$ |
| 0 | 2 | ... | $x_j-2$ | $x_j-1$ | $x_j+1$ | ... | $n-4$ | $n-2$ |
| 2 | ... | $x_j-2$ | $x_j-1$ | $x_j+1$ | ... | $n-4$ | $n-2$ | $n$ |
| 0 | ... | $x_j-4$ | $x_j-3$ | $x_j-1$ | ... | $n-6$ | $n-4$ | $n-2$ |
|   |   |   |   |   |   |   |   |   |
| 0 | 1 | 3 | ... | ... | ... | ... | ... | $n-2$ |

$j-1$ times (bracketing the middle rows)

Figure 7. The procedure of finding the prime form

This indicates that our assumption of $j \neq 1$ is false, and thus $j = 1$ and the first three elements in the prime form of $S_0^{(n,n_w),\cdot}$ are indeed $\{0,1,3,\dots\}$. Next, we show that this prime form violates Axiom II.

**Prime form and its violation of Axiom II**

Above, we analyzed pairs of continuous white keys, or the "tight part" in $S_0^{(n,n_w),prime}$. Next, we consider the "loose part." Consider the fragment of $S_0^{(n,n_w),prime}$ with configuration "01010." This fragment always exists; otherwise, the only possible key configuration of $S_0^{(n,n_w),prime}$ would be "001001…001," which would violate Axiom III by its periodic structure. This means that the prime form $S_0^{(n,n_w),prime}$ has the form $\{0,1,3,\dots,x_j-2,x_j,x_j+2,\dots,n-2\}$ for some $3 \leq j \leq n_w - 1$. We now claim that the key signature sequence $K_{n-2}^{(n,n_w),prime}$ contains $+2$ and thus violates Axiom II. To illustrate this, we add $n-2$, or subtract 2 from $S_0^{(n,n_w),prime}$. The elements "0" and "1" become $n-2$ and $n-1$, and we can conclude that $S_{n-2}^{(n,n_w),prime}$ is a version $S_0^{(n,n_w),prime}$ that has been shifted left twice and from which 2 has been subtracted.

Correspondingly, $S_{n-2}^{(n,n_w),prime} = \{1,\dots,x_j-4,x_j-2,x_j,\dots,n-4,n-2,n-1\}$ and thus $K_{n-2}^{(n,n_w),prime} = \{+1,\dots,\dots,\dots,+2,\dots,\dots,\dots,+1\}$. This means that the $(j-1)^{th}$ element of $K_{n-2}^{(n,n_w),prime}$ is $+2$, thus violating Axiom II.

# 5. PROPERTIES OF SOME KEY CONFIGURATIONS

**N=12**

For n = 12, there are three cyclically equivalent key configurations, namely, $S_0^{(12,7),1}$, $S_0^{(12,7),2}$, and $S_0^{(12,7),3}$. They are associated with piano scales starting with F, C, or G, respectively, or the "Lydian," "Ionian," and "Mixolydian" modes, respectively. In other words, the only viable keyboard layout for n=12 is the configuration that has been used for centuries. Another exciting result is that, among these solutions, the relative positions of the leading tones are identical. The ascending leading tone is always "B" (the white key following three black keys of equal distance), and the descending leading tone is always "F" (the white key before the three black keys of equal distance). The leading tones are thus called "shift invariant." However, the dominant and the subdominant are not shift invariant since they are defined in a fixed position of the scale.

An even more exciting property is related to the choice of a representative key configuration. Among equivalent configurations one can always find a configuration in which the ascending leading tone is the last note in the scale. Specifically, given $n, n_w, n_b$, the scale $S_0^{(n,n_w),2}$ is always the one that has the note with number $n-1$ as the ascending leading tone, as is the case in the usual 12-TET system. This is true for all possible $n, n_w, n_b$, thus making our axioms universal in all potential TET systems. Other interesting properties are as follows: (1) the scale $S_0^{(n,n_w),n_w-n_b+1}$ is always the one that has a note with number $1$ as the descending leading tone; (2) the note with number $n-2$ in the scale $S_0^{(n,n_w),1}$ is always the ascending leading tone; (3) the scale $S_0^{(n,n_w),1}$ is always the only one that has a black key as the last note; and (4) the scale $S_0^{(n,n_w),1}$ is always the inversion of the prime form. The position of the leading tone of $S_0^{(n,n_w),2}$ makes it the top keyboard-layout candidate for a given pair of $(n_w, n_b)$, and if not specified, $S_0^{(n,n_w)}$ refers to $S_0^{(n,n_w),2}$.

**N=19**

Perhaps the most studied equal-temperament system outside of 12-TET is 19-TET, which Joseph Yasser (1932) surveyed thoroughly, and for which he produced a keyboard layout (Figure 8).

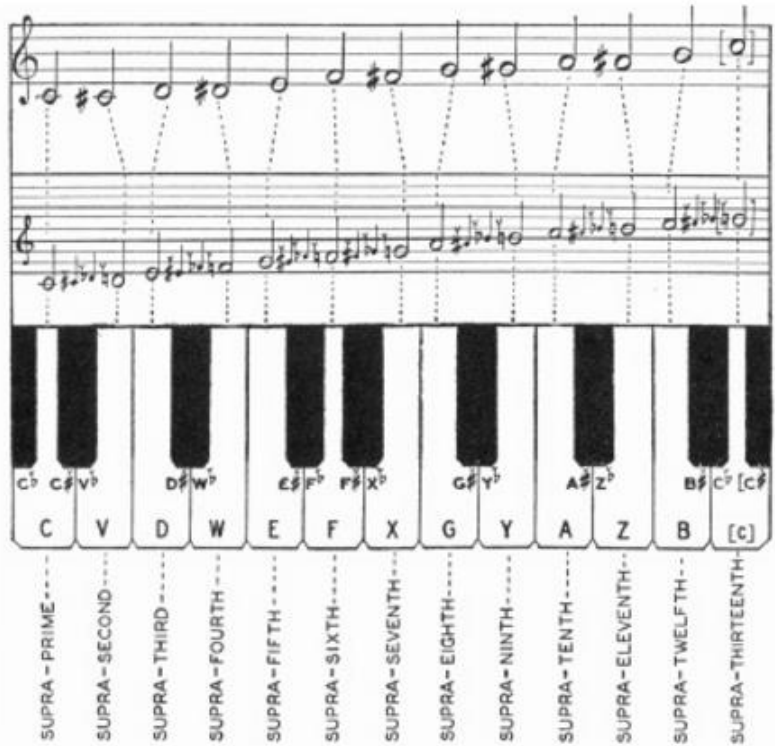

Figure 8. Joseph Yasser's 19 equal-temperament keyboard layout

We find that $S_0^{(19,12),1} = \{0,2,3,5,6,8,10,11,13,14,16,17\}$ is identical to Yasser's layout. Note that Yasser's starting point was evolving tones from 5+7 to 7+12 tones, and that he placed the seven black keys where the seven white keys are in 12-TET. Our exploration started out from a seemingly irrelevant perspective, but we have arrived at the same configuration as Yasser (1932).

Along with Yasser's keyboard-layout proposal, he also used the idea of the circle of fifths in 19-TET. The nearest interval to the perfect fifth in 19-TET is 11 semitones, so his circle of fifths was built with steps of 11 semitones. However, given our redefinition of the dominant in this article, we aim to see the circle built with steps that equal the distance between the tonic and the dominant to maintain consistency with a universal theory.

In this way, the mathematical properties of key signatures can be fully preserved. The dominant in 19-TET with $n_w = 12$ is the note with number 8, and the frequency ratio of the dominant to the tonic is $2^{8/19} = 1.339$, which is close to both the perfect fourth in equal temperament $(2^{5/12} = 1.335)$, and the perfect fourth in just intonation $(4/3 = 1.333)$. Compared to $2^{11/19} = 1.494$, which is 0.416% lower than $3/2$ and 0.304% lower than $2^{7/12}$, $2^{8/19} = 1.339$ is 0.418% higher than $4/3$ and 0.304% higher than $2^{5/12}$. With regards to the deviation, for tuning reasons there is no difference whether the step size is 8 or 11 semitones in 19-TET to finish the circle. However, a step size of $n_w^{-1} = 8$ is consistent with potentially universal rules. Given that 11 is the subdominant, choosing steps with interval 11 is equivalent to the reverse traversal of the circle, and thus Yasser's theory holds true. However, such a circle of key signature should not be called the circle of "fifths" because an interval of 8 is actually a sixth ({0,2,3,5,6,8}, or "supra-prime" to "supra-sixth", as Yasser defined). Therefore, it should be called the circle of sixths (in the sense of the diatonic scale in 19-TET), or the circle of perfect fourths (its frequency ratio approaches perfect fourths in 12-TET).

In addition, to be consistent with the convention in the 12-TET system, we suggest that the configuration be the one that has the ascending leading tone as the last note in the scale, that is, $S_0^{(19,12),2} = \{0,2,3,5,6,8,10,11,13,14,16,18\}$.

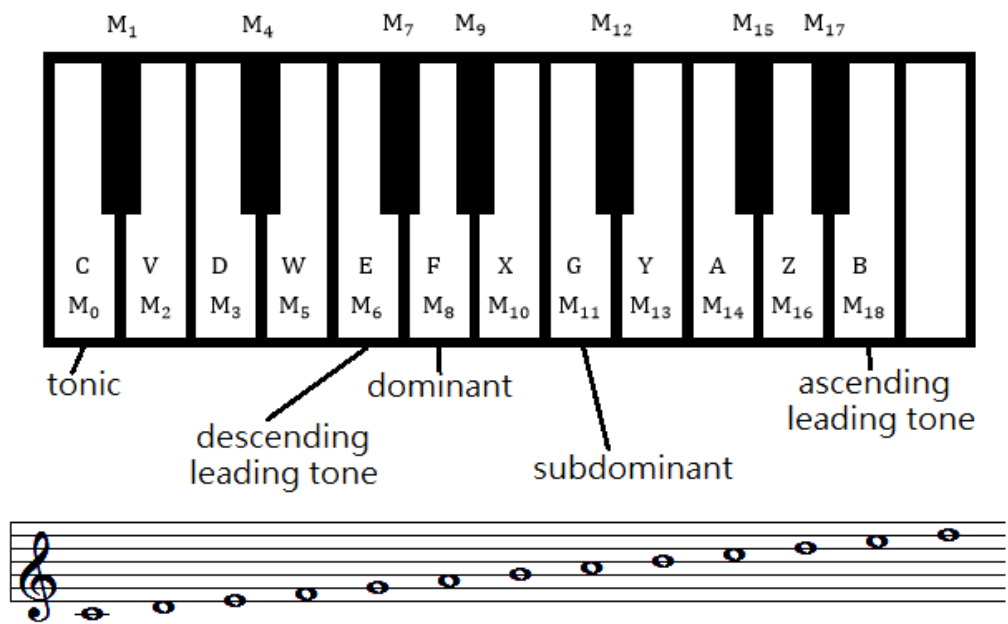

Figure 9. Our proposed 19 equal-temperament keyboard layout

This configuration differs from Yasser's in only one key, namely, position $M_{18}$ has been exchanged with $M_{17}$. Note that the nomenclature for notes (C, V, D, W, E,…, B) follows Yasser's naming system, so the note "B" differs by one semitone between the two configurations.

The staff that we adopt contains seven lines, chosen because it equals the number of black keys. By the first axiom of key arrangement, we have $n_b \geq \left\lceil \frac{1}{3}n \right\rceil$, so the interval of $n_b$ lines is roughly $2n_b$, which is greater than $n_w$. In other words, $n_b$ lines can contain at least an octave and hence our choice of staff seems reasonable.

The key signature matrix of scale $S_0^{(19,12),2}$ is

$$\mathbf{K}^{(19,12),2} = \begin{pmatrix} 0 & 0 & 1 & 0 & 1 & 0 & 0 & 1 & 0 & 0 & 1 & 0 & 1 & 0 & 0 & 1 & 0 & 0 & 1 \\ 0 & -1 & 0 & 0 & 1 & 0 & -1 & 0 & 0 & -1 & 0 & 0 & 1 & 0 & -1 & 0 & 0 & -1 & 0 \\ 0 & 0 & 1 & 0 & 1 & 1 & 0 & 1 & 0 & 0 & 1 & 0 & 1 & 1 & 0 & 1 & 0 & 0 & 1 \\ 0 & -1 & 0 & 0 & 1 & 0 & 0 & 1 & 0 & -1 & 0 & 0 & 1 & 0 & 0 & 1 & 0 & -1 & 0 \\ 0 & 0 & 1 & 0 & 1 & 1 & 0 & 1 & 1 & 0 & 1 & 0 & 1 & 1 & 0 & 1 & 1 & 0 & 1 \\ 0 & -1 & 0 & 0 & 1 & 0 & 0 & 1 & 0 & 0 & 1 & 0 & 1 & 0 & 0 & 1 & 0 & 0 & 1 \\ 0 & -1 & 0 & -1 & 0 & 0 & -1 & 0 & 0 & -1 & 0 & 0 & 1 & 0 & -1 & 0 & 0 & -1 & 0 \\ 0 & 0 & 1 & 0 & 1 & 0 & 0 & 1 & 0 & 0 & 1 & 0 & 1 & 1 & 0 & 1 & 0 & 0 & 1 \\ 0 & -1 & 0 & 0 & 1 & 0 & -1 & 0 & 0 & -1 & 0 & 0 & 1 & 0 & 0 & 1 & 0 & -1 & 0 \\ 0 & 0 & 1 & 0 & 1 & 1 & 0 & 1 & 0 & 0 & 1 & 0 & 1 & 1 & 0 & 1 & 1 & 0 & 1 \\ 0 & -1 & 0 & 0 & 1 & 0 & 0 & 1 & 0 & -1 & 0 & 0 & 1 & 0 & 0 & 1 & 0 & 0 & 1 \\ 0 & -1 & 0 & -1 & 0 & 0 & -1 & 0 & 0 & -1 & 0 & -1 & 0 & 0 & -1 & 0 & 0 & -1 & 0 \end{pmatrix}$$

And
$$\|\mathbf{K}^{(19,12),2}\| = (0 \ -7 \ 5 \ -2 \ 10 \ 3 \ -4 \ 8 \ 1 \ -6 \ 6 \ -1 \ 11 \ 4 \ -3 \ 9 \ 2 \ -5 \ 7).$$

Analogous to the circle of fifths in 12-TET, we provide the circle of sixths in 19-TET with scale $S_0^{(19,12),2} = \{0,2,3,5,6,8,10,11,13,14,16,18\}$ in Figure 10.

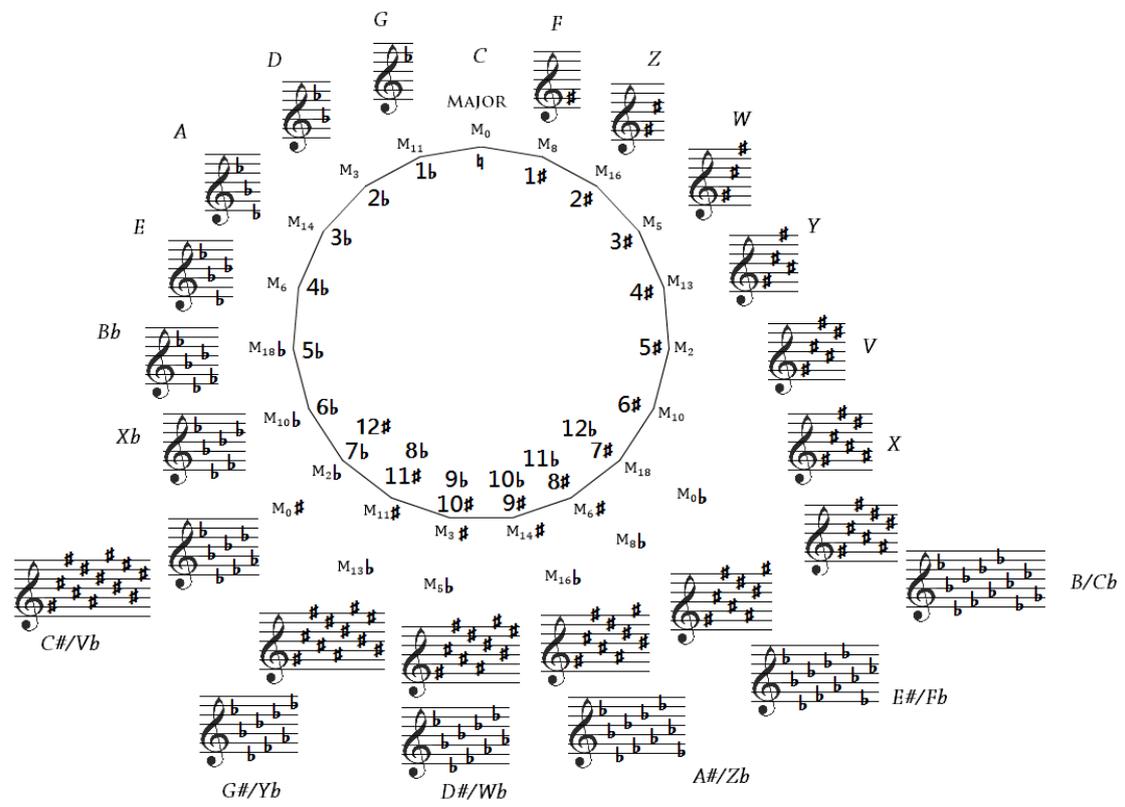

Figure 10. Circle of sixths in 19-TET

## 6. THE UNIVERSAL PERFECT FIFTH (FOURTH) CONSTANT

**Relations between the dominant and the perfect fifth**

So far only n = 12 and n = 19 have been discussed in part V; however, Table 6 contains further scales, such as n = 17, that remain unexplored. This is because they are not "elegant" enough. Firstly, we focus on the relation between the dominant and the perfect fifth. In 12-TET, we call the interval between the tonic and the dominant a "perfect" fifth because the ratio of their frequencies is $2^{7/12} = 1.498 \sim 3/2$, which is nearly a perfect ratio of two small integers. To determine whether the dominant really "dominates" the scale, we can examine this ratio in other equal-temperament systems.

Table 11. Frequency ratios of the dominant and the tonic in various scales

| Scale, $n, n_w, n_b$ | The dominant, $n_w^{-1}$ | Frequency ratio $2^{n_w^{-1}/n}$ |
|---|---|---|
| 12,7,5 | 7 | 1.498 |
| 13,8,5 | 5 | 1.306 |
| 13,7,6 | 2 | 1.113 |
| 14,9,5 | 11 | 1.724 |
| 15,8,7 | 2 | 1.097 |
| 16,9,7 | 9 | 1.477 |
| 17,11,6 | 14 | 1.770 |
| 17,10,7 | 12 | 1.631 |
| 17,9,8 | 2 | 1.085 |
| 18,11,7 | 5 | 1.212 |
| 19,12,7 | 8 | 1.339 |
| 19,11,8 | 7 | 1.291 |
| 19,10,9 | 2 | 1.076 |

As shown in Table 11, scales other than $S_0^{(12,7),2}$ and $S_0^{(19,12),2}$ usually have a frequency ratio of $2^{n_w^{-1}/n}$ that is not close to a perfect ratio.

**A miracle and Yasser's perspective of evolving tonality**

For a larger *n*, the table grows much faster, and the frequency ratio seems to become random and intractable, making it hard to find a good reason to keep track of specific pairs of $n_w, n_b$. Yasser, however, provided good guidance on finding interesting pairs of $n_w, n_b$.

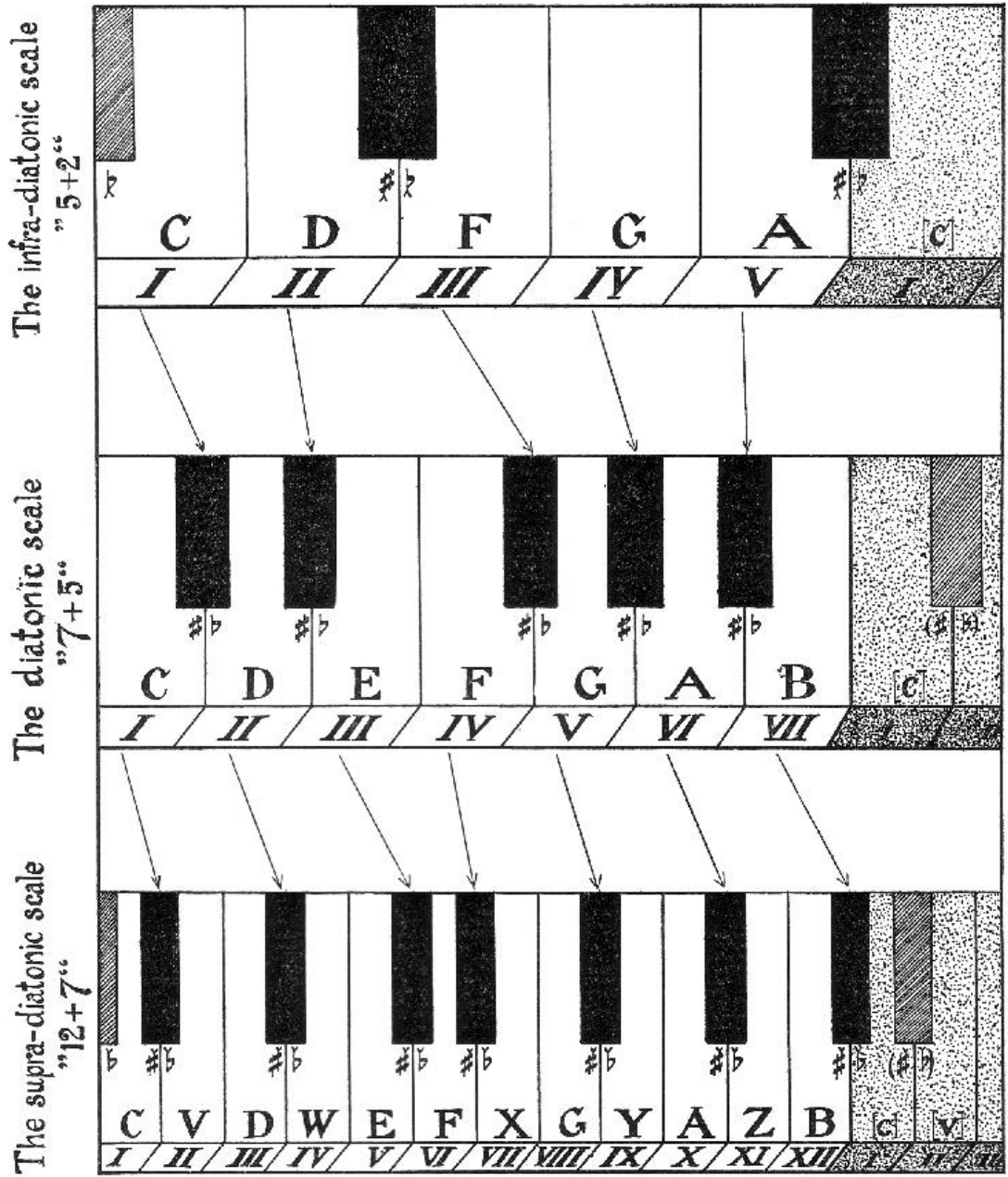

Figure 12. Yasser's evolving tonality

Figure 12 is reprinted from Carlo Serafini's review (2013) of Yasser's remarkable book "A theory of evolving tonality," and shows the process of evolving tonality. In short, the "white scale" of the deeper-level scales comes from the chromatic scale of the original scale, and the "white scale" of the original scale becomes the "black scale" of the deeper scale. This mechanism gives us a good criterion to select the next pair of $(n_w, n_b)$ with $(n_w', n_b') = (n_w + n_b, n_w)$. In Table 13, we list such interesting pairs of $(n_w, n_b)$.

Table 13. Dominants in evolved scales and the frequency ratios to their tonics

| Scale, $n, n_w, n_b$ | The dominant, $n_w^{-1}$ | Frequency ratio $2^{n_w^{-1}/n}$ |
|---|---|---|
| 7,5,2 | 3 | 1.346 |
| 12,7,5 | 7 | 1.498 |
| 19,12,7 | 8 | 1.339 |
| 31,19,12 | 18 | 1.496 |
| 50,31,19 | 21 | 1.338 |
| 81,50,31 | 47 | 1.495 |
| 131,81,50 | 55 | 1.338 |
| 212,131,81 | 123 | 1.495 |
| 343,212,131 | 144 | 1.338 |

It is interesting that $2^{n_w^{-1}/n}$ has two values, approximately 4/3 and 3/2, as the attractors. This is in fact highly surprising because the inverse element is highly unstable with merely a small change in $n_w$. For example, in $(n_w, n_b) = (81,50)$, $81^{-1} \equiv 55 \ (mod\ 131)$, but $80^{-1} \equiv 113 \ (mod\ 131)$ and $82^{-1} \equiv 8 \ (mod\ 131)$. The mapping from $n_w$ to $n_w^{-1}$ does not preserve ordering and vicinity. The tendency of $2^{n_w^{-1}/n}$ to attach to two fixed values is therefore not a coincidence and can be proved mathematically.

Firstly, $n_w$ is always coprime to $n$. Recall the Euclidean algorithm for computing the greatest common divisor (GCD) of two numbers. This algorithm is effective because of its nature to "divide and conquer," and the fact that $gcd(a, b) = gcd(b, a - b)$. Due to the method by which Yasser generated $(n_w, n_b)$ pairs, we have $gcd(n, n_w) = gcd(n_w, n - n_w) = gcd(n_w, n_b)$, and by repeating this process, we finally get $gcd(n, n_w) = gcd(5,2) = 1$. Therefore, $n_w$ is always coprime to $n$, so $n_w^{-1}$ always exists and is always unique.

It is somewhat confusing to write $n_w, n_b$ here because $n_b$ takes the place of $n_w$ after each turn of evolution. Therefore, we define a Fibonacci-like sequence $W_k$, with the property $W_{k+2} = W_{k+1} + W_k$ and the initial conditions that $W_1 = 2$ and $W_2 = 5$. If needed, $W_0$ can be defined as $W_0 = W_2 - W_1 = 3$. Now, every pair of $(W_{k+1}, W_k)$ is a pair of $(n, n_w)$, and finding $n_w^{-1}$ can be regarded as seeking $V_k = W_k^{-1} (mod\ W_{k+1})$.

Table 14. The relation between $W_k$, $V_k$, and the quotient $U_k$

| $W_k$ | $V_k = W_k^{-1} (mod\ W_{k+1})$ | $W_k V_k = U_k W_{k+1} + 1$ | $U_k$ |
|---|---|---|---|
| $W_1 = 2$ | $V_1 = 3$ | $W_1 V_1 = W_2 + 1$ | $U_1 = 1$ |
| $W_2 = 5$ | $V_2 = 3$ | $W_2 V_2 = 2W_3 + 1$ | $U_2 = 2$ |
| $W_3 = 7$ | $V_3 = 7$ | $W_3 V_3 = 4W_4 + 1$ | $U_3 = 4$ |
| $W_4 = 12$ | $V_4 = 8$ | $W_4 V_4 = 5W_5 + 1$ | $U_4 = 5$ |
| $W_5 = 19$ | $V_5 = 18$ | $W_5 V_5 = 11W_6 + 1$ | $U_5 = 11$ |
| $W_6 = 31$ | $V_6 = 21$ | $W_6 V_6 = 13W_7 + 1$ | $U_6 = 13$ |
| $W_7 = 50$ | $V_7 = 47$ | $W_7 V_7 = 29W_8 + 1$ | $U_7 = 29$ |
| $W_8 = 81$ | $V_8 = 55$ | $W_8 V_8 = 34W_9 + 1$ | $U_8 = 34$ |
| $W_9 = 131$ | $V_9 = 123$ | $W_9 V_9 = 76W_{10} + 1$ | $U_9 = 76$ |
| $W_{10} = 212$ | $V_{10} = 144$ | $W_{10} V_{10} = 89W_{11} + 1$ | $U_{10} = 89$ |

Proving directly that $2^{n_w-1}/n$ has two attractors is not easy. To simplify the proof, let us introduce an auxiliary number $U_k$. As long as we can find $0 < V_k, U_k < W_{k+1}$ such that $W_k V_k = U_k W_{k+1} + 1$ holds, then $V_k = W_k^{-1} (mod\ W_{k+1})$. Since $\gcd(W_k, W_{k+1}) = 1$, the existence and uniqueness of $V_k$ are guaranteed by elementary properties of rings of integers. In the following, we aim to find the general formula for $V_k$ by observing cases for some small $k$ in Table 14.

Let us consider cases separately for odd and even $k$. For odd $k$, $V_k$ can be obtained from $V_{k+2} = V_{k+1} + V_k$ with the initial conditions that $V_0 = 1$ and $V_1 = 3$, resulting in $V_1 = 3$, $V_2 = 4$, $V_3 = 7$, $V_4 = 11$, $V_5 = 18$, and so on. In the same way in which $U_k$ can be obtained with the initial conditions that $U_0 = 2$ and $U_1 = 1$ for odd $k$, $V_k$ can be obtained with the initial conditions of $V_0 = 1$ and $V_1 = 2$ for even $k$, and $U_k$ can be obtained with the initial conditions of $U_0 = 1$ and $U_1 = 1$ for even $k$.

To obtain the general formula for these Fibonacci-like sequences, we can write the relations of $W_{k+2}$, $W_{k+1}$, and $W_k$ as

$$\begin{bmatrix} W_{k+2} \\ W_{k+1} \end{bmatrix} = \begin{bmatrix} 1 & 1 \\ 1 & 0 \end{bmatrix} \begin{bmatrix} W_{k+1} \\ W_k \end{bmatrix} = A\vec{w_k}, \text{ with } \vec{w_0} = \begin{bmatrix} 2 \\ 3 \end{bmatrix}.$$

To find $W_k$ for arbitrary $k$, we have to compute $\vec{w_k} = A^k \vec{w_0}$. Note that $A$ has two linearly independent eigenvectors (non-normalized) $\vec{e_1} = \begin{bmatrix} \frac{1+\sqrt{5}}{2} \\ 1 \end{bmatrix}$ and $\vec{e_2} = \begin{bmatrix} \frac{1-\sqrt{5}}{2} \\ 1 \end{bmatrix}$, for $A\vec{e_1} = \frac{1+\sqrt{5}}{2} \vec{e_1} = \lambda_1 \vec{e_1}$ and $A\vec{e_2} = \frac{1-\sqrt{5}}{2} \vec{e_2} = \lambda_2 \vec{e_2}$.

Now we decompose $\overrightarrow{w_0}$ to $\overrightarrow{e_1}$ and $\overrightarrow{e_2}$, having $\overrightarrow{w_0} = \begin{bmatrix}2\\3\end{bmatrix} = \frac{15+\sqrt{5}}{10}\overrightarrow{e_1} + \frac{15-\sqrt{5}}{10}\overrightarrow{e_2}$; thus, $\overrightarrow{w_k} = A^k\overrightarrow{w_0} = \left(\frac{1+\sqrt{5}}{2}\right)^k \frac{15+\sqrt{5}}{10}\overrightarrow{e_1} + \left(\frac{1-\sqrt{5}}{2}\right)^k \frac{15-\sqrt{5}}{10}\overrightarrow{e_2}$, and we have

$$W_k = \left(\frac{1+\sqrt{5}}{2}\right)^k \frac{15+\sqrt{5}}{10} + \left(\frac{1-\sqrt{5}}{2}\right)^k \frac{15-\sqrt{5}}{10}.$$

By the same way, we have

$$V_{2m-1} = \left(\frac{1+\sqrt{5}}{2}\right)^{2m} + \left(\frac{1-\sqrt{5}}{2}\right)^{2m},$$

$$V_{2m} = \left(\frac{1+\sqrt{5}}{2}\right)^{2m} \frac{5+3\sqrt{5}}{10} + \left(\frac{1-\sqrt{5}}{2}\right)^{2m} \frac{5-3\sqrt{5}}{10},$$

$$U_{2m-1} = \left(\frac{1+\sqrt{5}}{2}\right)^{2m-1} + \left(\frac{1-\sqrt{5}}{2}\right)^{2m-1},$$

$$U_{2m} = \left(\frac{1+\sqrt{5}}{2}\right)^{2m} \frac{5+\sqrt{5}}{10} + \left(\frac{1-\sqrt{5}}{2}\right)^{2m} \frac{5-\sqrt{5}}{10}, \text{ for } m = 1,2,3,\ldots$$

Finally, we check the value of $W_k V_k - U_k W_{k+1}$:

For odd $k$:

$$\left[\left(\frac{1+\sqrt{5}}{2}\right)^{2m-1} \frac{15+\sqrt{5}}{10} + \left(\frac{1-\sqrt{5}}{2}\right)^{2m-1} \frac{15-\sqrt{5}}{10}\right]\left[\left(\frac{1+\sqrt{5}}{2}\right)^{2m} + \left(\frac{1-\sqrt{5}}{2}\right)^{2m}\right] -$$

$$\left[\left(\frac{1+\sqrt{5}}{2}\right)^{2m} \frac{15+\sqrt{5}}{10} + \left(\frac{1-\sqrt{5}}{2}\right)^{2m} \frac{15-\sqrt{5}}{10}\right]\left[\left(\frac{1+\sqrt{5}}{2}\right)^{2m-1} + \left(\frac{1-\sqrt{5}}{2}\right)^{2m-1}\right]$$

$$= \frac{15+\sqrt{5}}{10}\left(\frac{1+\sqrt{5}}{2}\right)^{2m-1}\left(\frac{1-\sqrt{5}}{2}\right)^{2m-1}\left(\frac{1-\sqrt{5}}{2} - \frac{1+\sqrt{5}}{2}\right) +$$

$$\frac{15-\sqrt{5}}{10}\left(\frac{1-\sqrt{5}}{2}\right)^{2m-1}\left(\frac{1+\sqrt{5}}{2}\right)^{2m-1}\left(\frac{1+\sqrt{5}}{2} - \frac{1-\sqrt{5}}{2}\right)$$

$$= \frac{15+\sqrt{5}}{10}(-1)(-\sqrt{5}) + \frac{15-\sqrt{5}}{10}(-1)(\sqrt{5}) = 1.$$

For even $k$:

$$\left[\left(\frac{1+\sqrt{5}}{2}\right)^{2m} \frac{15+\sqrt{5}}{10} + \left(\frac{1-\sqrt{5}}{2}\right)^{2m} \frac{15-\sqrt{5}}{10}\right]\left[\left(\frac{1+\sqrt{5}}{2}\right)^{2m} \frac{5+3\sqrt{5}}{10} + \left(\frac{1-\sqrt{5}}{2}\right)^{2m} \frac{5-3\sqrt{5}}{10}\right] -$$

$$\left[\left(\frac{1+\sqrt{5}}{2}\right)^{2m+1} \frac{15+\sqrt{5}}{10} + \left(\frac{1-\sqrt{5}}{2}\right)^{2m+1} \frac{15-\sqrt{5}}{10}\right]\left[\left(\frac{1+\sqrt{5}}{2}\right)^{2m} \frac{5+\sqrt{5}}{10} + \left(\frac{1-\sqrt{5}}{2}\right)^{2m} \frac{5-\sqrt{5}}{10}\right]$$

$$= \left(\frac{1+\sqrt{5}}{2}\right)^{4m}\left[\frac{15+\sqrt{5}}{10}\frac{5+3\sqrt{5}}{10} - \frac{15+\sqrt{5}}{10}\frac{5+\sqrt{5}}{10}\frac{1+\sqrt{5}}{2}\right] + \left(\frac{1-\sqrt{5}}{2}\right)^{4m}\left[\frac{15-\sqrt{5}}{10}\frac{5-3\sqrt{5}}{10} - \right.$$

$$\left.\frac{15-\sqrt{5}}{10}\frac{5-\sqrt{5}}{10}\frac{1-\sqrt{5}}{2}\right] + \left(\frac{1+\sqrt{5}}{2}\right)^{2m}\left(\frac{1-\sqrt{5}}{2}\right)^{2m}\left[\frac{15+\sqrt{5}}{10}\frac{5-3\sqrt{5}}{10} + \frac{15-\sqrt{5}}{10}\frac{5+3\sqrt{5}}{10} - \right.$$

$$\left.\frac{15+\sqrt{5}}{10}\frac{5-\sqrt{5}}{10}\frac{1+\sqrt{5}}{2} - \frac{15-\sqrt{5}}{10}\frac{5+\sqrt{5}}{10}\frac{1-\sqrt{5}}{2}\right]$$

- 29 -

$$= \left(\tfrac{1+\sqrt{5}}{2}\right)^{4m}(0) + \left(\tfrac{1-\sqrt{5}}{2}\right)^{4m}(0) + (-1)^{2m}\left[\tfrac{120}{100} - \tfrac{40}{200}\right] = 1.$$

Therefore, $W_k V_k = U_k W_{k+1} + 1$ is true, and $V_k = W_k^{-1} \pmod{W_{k+1}}$ has the general formula

$$V_{2m-1} = \left(\tfrac{1+\sqrt{5}}{2}\right)^{2m} + \left(\tfrac{1-\sqrt{5}}{2}\right)^{2m},$$

$$V_{2m} = \left(\tfrac{1+\sqrt{5}}{2}\right)^{2m}\tfrac{5+3\sqrt{5}}{10} + \left(\tfrac{1-\sqrt{5}}{2}\right)^{2m}\tfrac{5-3\sqrt{5}}{10}, \text{ for } m = 1,2,3,\ldots$$

The two attractors of $2^{n_w{}^{-1}/n}$ can be written as $\lim\limits_{m\to\infty} 2^{\frac{V_{2m-1}}{W_{2m}}}$ and $\lim\limits_{m\to\infty} 2^{\frac{V_{2m}}{W_{2m+1}}}$, and their values are

$$\lim_{m\to\infty} 2^{\frac{V_{2m-1}}{W_{2m}}} = 2^{\lim_{m\to\infty}\frac{V_{2m-1}}{W_{2m}}} = 2^{\frac{1}{\frac{15+\sqrt{5}}{10}}} = 2^{\frac{15-\sqrt{5}}{22}} = 1.49503444953\ldots$$

$$\lim_{m\to\infty} 2^{\frac{V_{2m}}{W_{2m+1}}} = 2^{\frac{\frac{5+3\sqrt{5}}{10}}{\frac{15+\sqrt{5}}{10}\cdot\frac{1+\sqrt{5}}{2}}} = 2^{\frac{7+\sqrt{5}}{22}} = 1.33776181588\ldots$$

Because the product of the two constants is exactly 2, which resembles the cases of the perfect fifth and fourth in 12-TET, respectively, we call them the "universal perfect fifth constant" and the "universal perfect fourth constant," respectively. We conclude that the interval between the dominant and the tonic is almost always either a perfect fifth or a perfect fourth, as long as the scale is "good," and Yasser always produced good scales!

# 7. CONCLUSION

To date, the mathematical properties of scales have been thoroughly studied in diatonic-set theory. Instead of pointing out the new properties of these well-studied scales, this article offers viewpoints from the scoring system. Starting with observations of regular piano- and musical-notation systems, we proposed several axioms, found the correct 12-TET solutions, and generalized studied aspects to general equal-temperament systems. While this study does not focus on specific proofs, we logically develop our considerations and attempt to tell a story that makes sense for both mathematicians and musicians, explaining the logic behind musical conventions. The main point of this article is to unify and to suggest scale configurations, degrees, and ideas such as the circle of fifth in general TET systems that are consistent to fundamental axioms. This article ends with a beautiful constant, which can be ascribed to remarkable ideas from brilliant researchers over time.